\documentclass{elsart}

\usepackage{graphicx}

\usepackage{amsfonts}  

\usepackage{amsmath}

\def\R{$\mathbb{R}$}

\def\F{$\mathbb{F}$}
\def\GF{$\mathbb{GF}$}

\newcommand{\mqed}{\hfill \ensuremath{\Box}}

\newcounter{mycount}
\newcounter{mydcount}
\newcounter{mytcount}
\newcounter{mylcount}
\newcounter{myccount}

\newcommand\mydef[1]{%
  \stepcounter{mydcount}
   \par\noindent {\bf Definition \themydcount: {\it #1}}
}

\newcommand\mylem[1]{%
 \ \stepcounter{mylcount}
   \par\noindent {\bf Lemma \themylcount: {\it #1}}
}

\newcommand\mythm[1]{%
 \ \stepcounter{mytcount}
   \par\noindent {\bf Theorem \themytcount: {\it #1}}
}

\newcommand\myproof{%
\ 
   \par\noindent {\it Proof: }
 \ 
}

\newcommand\myprob[1]{%
\ \stepcounter{mycount}
   \par\noindent {\bf Problem \themycount: #1}
 \\ \\
}

\newcommand\mycor[1]{%
\ \stepcounter{myccount}
   \par\noindent {\bf Corollary \themyccount: {\it #1}}
}

\begin{document}

\begin{frontmatter}
%\begin{large}

\title{Flexibility of Bricard's linkages and other structures via resultants and computer algebra}
\author{Robert H. Lewis} \\
\address{Fordham University, New York, NY 10458, USA \\
http://fordham.academia.edu/RobertLewis}
\author{Evangelos A. Coutsias} \\
\address{Stony Brook University, Stony Brook, NY 11794, USA \\
http://www.ams.stonybrook.edu/$\sim$coutsias}

\maketitle

\begin{abstract}
Flexibility of structures is extremely important for chemistry and robotics.  Following our earlier work, we study flexibility using polynomial equations, resultants, and a symbolic algorithm of our creation that analyzes the resultant.  We show that the software solves a classic arrangement of quadrilaterals in the plane due to Bricard.  We fill in several gaps in Bricard's work and discover new flexible arrangements that he was apparently unaware of.   This provides strong evidence for the maturity of the software, and is a wonderful example of mathematical discovery via computer assisted experiment.
\end{abstract}

\end{frontmatter}
\section{Introduction}

This project results from the convergence of four topics: systems of polynomial equations, flexibility of two  and three 
dimensional objects, computational chemistry,  and computer algebra.  It also has application to robotics \cite{man}, \cite{hust}.

We have developed software to detect flexibility in certain structures that are generically rigid.  It is based on symbolic computation of polynomials and rational functions, not numerical computing.  We previously reported on earlier stages of this research in \cite{adg} and \cite{fox}.  Since then, the software has been enormously improved in both power and efficiency, to the point where it not only discovers the previously known modes of flexibility of a classic structure due to R. Bricard, but discovers new modes apparently unknown to him.  

We are mostly concerned with the framework in Figure 1.  It is a system of seven bars, joined at the nine junctions shown by rotational joints, allowing free rotation in the plane.  It is {\it generically rigid.}  This follows from a general theorem in kinematics \cite{hunt} by which the mobility
(number of degrees of freedom of relative motion) of a linkage system
is given by the relation
\[ {\cal{M}} = 3(n-g-1) + \sum_{i=1}^g f_i \]
where $n$ is the number of members, $g$ is the number of joints and $f_i$ is
the mobility at joint $i$.
For the system in Figure 1,
comprised entirely of rigid rods with rotatable joints (with $n=7$, $g=9$ and
$f_1,...,f_9=1$)
this gives ${\cal{M}} = 0$.

When $M > 0$ the system is generically flexible.  When $M = 0$ it is generically rigid or {\it determined}.  We wish to discover cases, by means of particular relations existing between its edges, that
determinacy (rigidity) ceases to hold. Then the framework will be deformable (flexible).

Flexibility is an intuitive concept.  Imagine a triangle made of three stiff rods joined with movable hinges. The formula above confirms the clear intuition that the structure is obviously rigid.  In the same way, a quadrilateral in the plane is obviously flexible, and for it $M > 0$. ($M$ = 1).  

   In computational chemistry,  protein folding has been a major research topic for a number of years \cite{thor}, \cite{cout}, \cite{cout2}.
Molecules can fold because they are flexible.  Simple examples are easily built from a few plastic balls and rods.  
In 1812, Cauchy considered flexibility of three dimensional polyhedra with triangular facets (similar to a geodesic
dome)  where each joint can pivot or hinge.  He proved that if the polyhedron is convex  it cannot be flexible; it must
be rigid \cite{cauch}.   In 1896 Bricard \cite{bric} tried to find non-convex flexible polyhedra by looking at one of the smallest possible cases, octahedra.  He partially succeeded: his flexible octahedra are not imbedible in actual three-space because some of their facets intercross.  He also described the system of three quadrilaterals in the plane (Figure 1) whose motion is algebraically equivalent to the octahedra.  

People came to believe that there were no flexible polyhedra at all.  But in 1978 Robert Connelly, building on Bricard, 
astonished them by finding a non-convex one \cite{con}, and soon models appeared of a simpler flexible structure
\cite{stef}, \cite{mak}. 

Our approach is to describe the geometry of the object or molecule with a set of
multivariate polynomial equations.  Solving a system of multivariate polynomial equations is a
classic, difficult problem.  The approach via resultants was pioneered by Bezout \cite{bik}, Dixon \cite{buse}, \cite{dix}, \cite{kap},  Sylvester \cite{cox}, 
and others.  The resultant appears as a factor of the determinant  of a matrix containing multivariate
polynomials.  Computing it can be quite a challenge \cite{ls}, but we developed methods to do so \cite{lh}.  Once we have the resultant, we described \cite{adg} an algorithm  we call {\it Solve} that examines the resultant and determines ways that the structure can be flexible.

\vspace*{1mm}
\begin{figure} 
\centering
\includegraphics[scale=0.32]{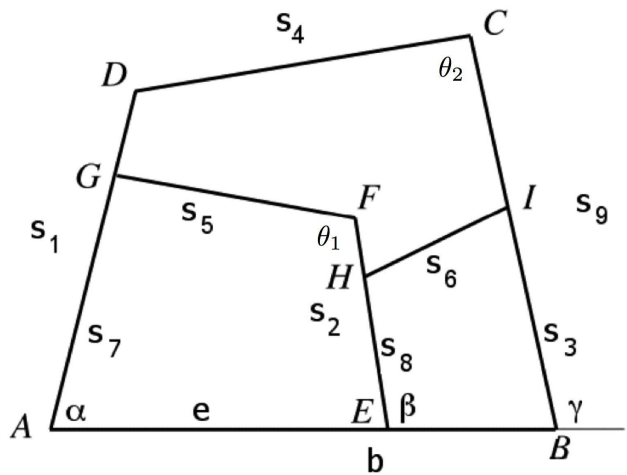}
\caption[ ]{Bricard's quadrilaterals, showing labeled sides and base angles.} 
\end{figure}
\vspace*{-2mm}

We discovered in this way some of the conditions of flexibility
for Bricard's arrangement of quadrilaterals in \cite{bric} that is algebraically equivalent to the octahedra. {\it Solve} was greatly improved by Fox \cite{fox} and more recently again by us.    It is at least 500 times faster on Bricard's quadrilaterals than in 2008, and now finds all three of Bricard's ways the quadrilaterals can flex (this was not true in \cite{adg}).   Surprisingly, {\it Solve} has discovered new flexible arrangements of the quadrilaterals that were apparently not anticipated by Bricard.

The main point of this paper is that our new algebraic and algorithmic solution of Bricard's quadrilaterals demonstrates that the software has matured to the point where one can confidently use it on more complex structures, such as molecules. 

\section{Basic setup and strategy}

All computations here were done with Lewis's computer algebra system  {\it Fermat} \cite{lfer},
which excels as polynomial and matrix computations \cite{rob}.

   As explained in the introduction, we are primarily concerned here with the analysis of the flexibility of a certain structure of Bricard consisting of seven rigid rods forming three quadrilaterals in the plane (Figure 1).  We need to establish that our software can find all the flexible cases.  It order to show why we are confident in this, we will present proofs paralleling some of those of Bricard \cite{bric}, but ours are quite different and more algebraic.

 Other than historical motivation, why should we concentrate on this arrangement of rods?  As remarked in the introduction, a quadrilateral in the plane $AD, DC, CB, BA$ is obviously flexible ($M=1$; see Figure 1).  Imagine that $AB$ is fixed. As $AD$ pivots about $A$, the angles $\alpha$ and $\gamma$ take on a continuum of values.  If we add $GF$ and $FE$ we have two nested quadrilaterals, and the structure remains obviously flexible ($M=1$), and $\beta$ also takes on a continuum of values, unless we set $F = A$, so $GF=GA$ and $FE=AE, s_5 = s_7, s_2 = e$; then $\beta$ would be constantly $\pi$ during the flex.  That is a {\it degenerate} case.  The addition of one more rod or ``brace" $HI$  produces a third quadrilateral $EHIB$. The structure is now generically rigid ($M=0$).  However, it can be made flexible in several ways.  A degenerate way to do so is to simply place $HI$ on top of $EB$, so $HI = EB, s_3 = s_8 = 0, s_6 = b-e$.  We are not concerned with such degeneracies here.\footnote{Informally, degeneracy means a side is 0, or an angle is constant during the flex.}   Far more interesting is to choose the lengths of the rods (sides) so that each quadrilateral is a parallelogram.  Obviously,  the system is then flexible. This is one of the cases we analyze below (section four).

This is our goal:  {\it non-generic flexibility.}  The system of Figure 1 is one of the simplest to examine for non-degenerate flexibility, and was thought by Bricard to be ``equivalent" to two octahedra in three dimensions.

 Our strategy is to describe the arrangement by a system of three polynomial equations, where the parameters are the lengths of the sides, and the three variables  represent the tangents of certain 
angles\footnote{Technically, the tangent of one-half the angle.} in the structure.   Using resultants, we show that flexibility implies that each of the three tangents is a rational function of the other two, and analyze when one tangent can be a rational function of only one other.  We thereby derive the three flexible cases that were defined by Bricard, but with new subcases.

% \noindent
 As in Bricard's paper on flexible octahedra, elementary geometry and trigonometry lead to a system of three polynomial equations in three variables $t_1, t_2, t_3$ and fifteen parameters (details in next paragraph), namely
\begin{equation} \label{eq:e1}
   a_1 t_1^{\, 2} t_2^{\, 2} + d_1 t_2^{\, 2} + 2 c_1 t_1 t_2 + b_1 t_1^{\, 2} + e_1 = 0 
\end{equation}
\vspace{-9mm}
\begin{equation}\label{eq:e2}
a_2 t_2^{\, 2} t_3^{\, 2} + d_2 t_3^{\, 2} + 2 c_2 t_2 t_3 + b_2 t_2^{\, 2} + e_2 = 0 
\end{equation}
\vspace{-6mm}
\begin{equation} \label{eq:e3} 
  a_3 t_1^{\, 2} t_3^{\, 2} + d_3 t_3^{\, 2} + 2 c_3 t_1 t_3 + b_3 t_1^{\, 2} + e_3 = 0  
\end{equation}

 The fifteen parameters are themselves simple polynomial functions of the sides of the flexing quadrilaterals, in such a way that $c_1 c_2 c_3 \ne 0$.  The other parameters might be 0.  The $t_i$ are half-angle tangents of angles $\alpha, \beta, \gamma$ in the quadrilaterals; see Figure 1.   (Cotangents could also be used, which has the effect of replacing $t_i$ with $t_i^{-1}$.  This will be used in Theorem 3.) The seven rigid rods are $AD, DC, CB, AB, GF, FE, HI$. The joints allow each rod to pivot freely in the plane. $AB$ remains fixed on the $x-$axis during pivoting, with $A$ at the origin.  We allow negative values for $s_3, s_8, s_2$, or $s_7$, so points $G, F, H, I$ might be below the $x$-axis. Angles $\theta_1$ and $\theta_2$ will be discussed later in Theorem 4. 
 
 The equations arise from Figure 1 using basic geometry and trigonometry.   For example, the coordinates of the point $D$ are $(s_1 \cos(\alpha), s_1$  $\sin(\alpha))$.  $C$ is $(b + s_9 \cos(\gamma), s_9 \sin(\gamma))$.  Therefore   
\begin{align*} s_4^2 = (b + s_9 \cos(\gamma) - s_1 \cos(\alpha))^2 + (s_9 \sin(\gamma) - s_1 \sin(\alpha))^2
\end{align*}  
 One obtains three equations of this kind (using also $s_5$ and $s_6$) and three obvious equations of type $\sin^2(x) + \cos^2(x) = 1$.  Then use the well-known half angle tangent substitutions   
\begin{align*}
\sin(\alpha) = 2 t_1 /(1 + t_1^2)  \\
  \cos(\alpha) =(1 - t_1^2) /(1 + t_1^2) 
\end{align*}
 and so on (with $\beta, \gamma; t_2, t_3$) to form the three equations $(1)-(3)$.  The fifteen parameters become

\noindent
\vspace{-13mm}
%\addtocounter {equation} {1}
%\begin{align}
\noindent
\begin{eqnarray}  \label{eq:e4}
 a_1 =  (-s_2+e-s_5+s_7) (-s_2+e+s_5+s_7) \nonumber \\
b_1  =  (s_2 + e - s_5 + s_7) (s_2 + e + s_5 + s_7)    \nonumber  \\
c_1 = -4s_2\, s_7   & \ &   \nonumber    \\
d_1 = (-s_2 + e + s_5 - s_7) (-s_2 + e - s_5 - s_7)  & \ &   \nonumber \\
e_1 = (s_2 + e - s_7 - s_5) (s_2 - s_7 + s_5 + e) \nonumber \\\
a_2 = (-b + e + s_3 - s_8 - s_6) (-b + e + s_3 - s_8 + s_6)  \nonumber \\\
b_2 = (-b + e - s_3 + s_6 - s_8) (-b + e - s_3 - s_6 - s_8)  \nonumber \\\
c_2 = -4s_3\, s_8  & \ &   \nonumber \ \\
d_2 = (-b + e + s_8 + s_3 - s_6) (-b + e + s_8 + s_3 + s_6)   & \ &  \nonumber \\\
e_2 = (-b + e - s_3 + s_8 - s_6) (-b + e - s_3 + s_8 + s_6)  & \ &  \nonumber \ \\
a_3 = (b + s_1 - s_9 - s_4) (b + s_1 - s_9 + s_4)  \nonumber \\\
b_3 = (b + s_1 + s_4 + s_9) (b + s_1 - s_4 + s_9)   & \ &  \nonumber \  \\ 
c_3 = -4s_9\, s_1  & \ &   \nonumber \\\
d_3 = (-b + s_1 + s_9 + s_4) (-b + s_1 + s_9 - s_4)   & \ &  \nonumber \\
e_3 = (-b + s_1 - s_9 - s_4) (-b + s_1 - s_9 + s_4)  & \ &   \\  \nonumber \
\end{eqnarray}  \nonumber \

\vspace{-9mm}
\noindent
None of the sides $s_i, e, b$ is 0.  $b \ne e, s_7 \ne s_1, s_3 \ne s_9, s_2 \ne s_8$.  For convenience, we also define $s_{10} \equiv e, s_{11} \equiv b$, and we also refer to $s_1 - s_7$, $s_9 - s_3$, and $s_2 - s_8$ as ``sides".

As we discussed above, the  arrangement of quadrilaterals in Figure 1 is generically rigid.  That means, in spite of the flexible joints, if numerical values were assigned arbitrarily for the eleven sides, the angles $\alpha, \beta, \gamma$ would be uniquely determined.\footnote{At least up to sign or supplement. Some assignments would be impossible.}  The main task of this paper is 
%  \medskip
\myprob{Find conditions  on the sides under which the quadrilateral arrangement becomes flexible.}
Flexibility is marked mathematically by the three angles, and their half-angle tangents $t_i$, each taking on uncountably many values.  As remarked above, if all three quadrilaterals are parallelograms, the arrangement is flexible.  This means that under the substitutions $s_9 = s_1, b=s_4, s_2=s_7, s_5=e,s_8=s_3, s_6=b-e$,  not only are there common roots to the system of equations $(1)-(3)$, but there is a continuum of common roots; each $t_i$ is a never-constant continuous function $t_i: I \rightarrow$ \R. {\it Never-constant} means there is no open interval over which $t_i$ is constant.  Allowing that would create degenerate cases, which we do not discuss here.  (Some are discussed in \cite{adg}.)
We therefore have secondarily: 

\par\noindent {\bf Problem 1$'$: Find all conditions on the sides under which the quadrilateral arrangement becomes non-degenerate flexible.}
% \medskip

  To understand flexibility, we follow Bricard and ask    
\myprob{When is one of these variables, $t_2$, say, \\
%\begin{itemize}
%\item 
(1) a rational function of another $t_j$, or\\
%\item  
(2) a rational function of both of the other ones $t_1, t_3$?
% \end{itemize}
}
Using resultants, we will show in our Main Theorem (section 5) that flexibility always implies the second case.  The first case is referred to as {\it splitting.}\footnote{One may fairly ask for the motivation for Problem 2.  Recall that a similar question about the roots of a polynomial is the basic idea in Abel's analysis of the unsolvability of the quintic.}

   To make sense of ``rational function" we must discuss the ground field, \GF.  Let \F \ be a field. In many of our algebraic results, \F \ could be any field of characteristic not 2. However, eventually we will evaluate expressions like those in (\ref{eq:e4}) by substituting each parameter with an element of \F.  Therefore \F \ = a subfield of \R \ is appropriate.  We do not allow the sides to be arbitrary complex numbers.

   Then given \F,  we may first think of the ground ring as  \F$[a_1, \ldots, e_3]$ and the ground field  \GF \ as \F$(a_1, \ldots, e_3)$, the field of rational functions over  \F \ of the fifteen parameters.   However, the ground field is really \GF \ = \F$(s_1, s_2, \ldots, s_9, e, b)$, where the fifteen parameters are replaced with their definitions in (\ref{eq:e4}) above.  This means that large polynomials in the fifteen $a_1, \ldots, e_3$ must sometimes be thought of as even larger polynomials in the eleven 
$s_1, s_2, \ldots, s_9, e, b$.

  The evaluation homomorphism, obtained by substituting parameters with values in \F, can be thought of as a map from \GF \ to itself.  Also, when we speak of finding a solution to the system $(1)-(3)$, we understand as usual that the common root may lie in an extension field of the ground field, for example, a radical extension.

   We can now specify what we mean in Problem 1 by ``find conditions on the sides under which the quadrilateral arrangement becomes flexible."  We mean find substitutions of the form $s_i = p(s_1, s_2, \ldots, \hat{s_i}, \ldots)$, where $p \in$ \GF, so that the $t_i$ are continuous functions from some interval to \R.  We will show that this notion of ``condition" does indeed lead to both old and new flexibile arrangements.  That in turn suggests:

{\bf Problem 3:  Can all flexible cases be represented by a table of substitutions in this sense?}  

We will see in the conclusion that, very surprisingly, this is false, and we conjecture a modification of it.

 The rest of the paper is organized as follows: In section three we develop three lemmas to identify when an equation splits.  The various split cases are summarized in section four.  In section five we present the main theorem, which solves Problem 2 and says that if no equation splits, then every $t_i$ is a rational function of the other two.  In section six we complete the theory of the non-split case.  In section seven we describe the software results and two surprising new flexible cases for the quadrilaterals that were apparently unknown to Bricard.

In all of the following definitions, lemmas, and theorems we assume flexibility.  Some of them are true without this assumption, but we are not concerned with that.
\section{Splitting lemmas}  
\mydef{We say that one of the equations $(\ref{eq:e1})-(\ref{eq:e3})$ {\it splits} or {\it decomposes} if one of the $t_i$ in it can be expressed as a rational function of the other one.}

For simplicity, let's concentrate  on solving for $t_2$ in equation $(1)$.  Suppose in $(1)$ we have $a_1 = d_1 = 0$.  Then (1) reduces to 
\begin{equation} \label{eq: e5}
2 c_1 t_1 t_2 + b_1 t_1^{\, 2} + e_1 = 0 
\end{equation}
Since $c_1$ cannot be 0, we can solve this for $t_2$ and obtain a rational function, so $(1)$ would split.  This example is an important case in the following lemmas.

Assuming that $a_1 \ne 0$ or $d_1 \ne 0$, it is natural to solve for $t_2$ using the quadratic formula.  We have

\begin{equation}   \label{eq:e6}
t_2 = { -c_1 t_1 \pm \sqrt{c_1^2 t_1^2 - (a_1 t_1^2 + d_1)(b_1 t_1^2 + e_1)} \over a_1 t_1^2 + d_1 }
\end{equation}
\mydef{The polynomial under the square root sign in (\ref{eq:e6}) is called $F(t_1)$.}
\mylem{$t_2$ is a rational function of $t_1$ if and only if $F(t_1)$ is a perfect square in \F$[s_1, s_2, \ldots, s_9, e, b] \, [t_1]$.}
\myproof  If $a_1 = 0$ and $d_1 = 0$ the result is immediate.  So assume $a_1 \ne 0$ or $d_1 \ne 0$.

The ``if" part of the statement is obvious.  To prove the converse, assume that there is rational function $t_2 = f/g$  with $f, g \in$ \F$[s_1, s_2, \ldots, s_9, e, b] \, [t_1]$.  Inserting this into (\ref{eq:e6}) and clearing denominators yields 
$$  f \cdot (a_1 t_1^2 + d_1) = g \cdot ( -c_1 t_1 \pm \sqrt{F(t_1)} \,)   $$
Multiply it out, collect terms, and solve for $F(t_1)$:
$$  F(t_1) = p^2 / q^2 $$
where $p$ and $q$ are polynomials in \F$[s_1, s_2, \ldots, s_9, e, b] \, [t_1]$.  But that ring is a UFD.  By a standard argument with irreducible polynomials, $q^2$ divides $p^2$, so we are done.  \mqed
\mylem{If $F(t_1)$ is a perfect square in \F$[s_1, s_2, \ldots, s_9, e, b] \, [t_1]$ then $a_1b_1 = 0$ and $d_1 e_1 = 0$.}
\myproof 
Note that 
$$ F(t_1) = - a_1 b_1 t_1^4 + (c_1^2-a_1 e_1 - b_1d_1) t_1^2 - d_1 e_1$$ is
a polynomial in $t_1^2$.  If this is truly a quadratic in $t_1^2$ and a perfect square then its discriminant must be 0. But when the parameters $a_1, b_1, \ldots, e_1$ are replaced with their expressions (\ref{eq:e4}) in terms of the eleven sides $s_1, s_2, \ldots, e$, the discriminant simplifies enormously to $256 e^2 s_2^2 s_5^2 s_7^2$.\footnote{We recommend a computer algebra system for this computation.  However, Bricard did not have one!}  Therefore, the discriminant cannot be 0.  The only solution is that $F(t_1)$ is linear in $t_1^2$. Thus,  $a_1 =0$ or $b_1 = 0$.  

We now have that 
$$ F(t_1) =  (c_1^2-a_1 e_1 - b_1d_1) t_1^2 - d_1 e_1$$   If this is a perfect square, then it equals some $(At_1 + B)^2$.  As there is no linear term in $F(t_1)$, it must be that $A=0$ or $B=0$.  

We will prove that $A=0$ leads to a contradiction. We showed above that there are two cases to consider, $a_1=0$ or $b_1=0$.  Assume that $a_1=0$.  Then $A=0$  implies that $c_1^2=b_1d_1$.  Since $a_1 = (-s_2+e-s_5+s_7) (-s_2+e+s_5+s_7) = 0$ this in turn leads to two cases.  If $-s_2+e-s_5+s_7 = 0$, then $s_2+s_5=e+s_7$.  When this is plugged into the definitions of $c_1, b_1$ and $d_1$, we see after some computation that $c_1^2=b_1d_1$ reduces to $s_5+s_2=s_7$. But then $e=0$, impossible.  The second alternative, $-s_2+e+s_5+s_7=0$, leads in the same way to the same contradiction that $e=0$.
  
If $b_1=0$, the argument is analogous.  Now $A=0$  implies that $c_1^2=a_1 e_1$. $b_1 = (s_2 + e - s_5 + s_7) (s_2 + e + s_5 + s_7)=0$ so we again have two subcases.  Each leads to the contradiction $e=0$.

 Therefore, $A=0$ is impossible so it must be that $B=0$.  Therefore $d_1 e_1 = 0$.   
$\mqed$ 
\mylem{In equation $(i)$, $i = 1, 2, 3$, consider the six ways to choose a pair of $\{ a_i, b_i, d_i, e_i \}$, the four parameters that might be 0.  In all six cases, if that pair of parameters is 0, the equation splits.} 
\myproof We illustrate with $i=1$. The case of $ \{ a_1, d_1 \}$ was shown in the above example (\ref{eq: e5}) with $t_2$ a rational function of $t_1$. $ \{ a_1, b_1 \}$ is analogous, solving for $t_1$.  

For  $ \{ b_1, e_1 \}$  we have
$$  a_1 t_1^{\, 2} t_2^{\, 2} + d_1 t_2^{\, 2} + 2 c_1 t_1 t_2 = 0 $$
Since $t_2$ is a function taking on a range of values, it may be divided out and we obtain
$$  a_1 t_1^{\, 2} t_2 + d_1 t_2 + 2 c_1 t_1  = 0 $$
whence  we may solve for $t_2$ as 
$$ t_2 = {-2 \,c_1 t_1 \over a_1 t_1^2 + d_1} $$
This is valid unless both $d_1$ and $a_1$ are 0.  But if that were true, we would have $2 c_1 t_1 t_2 = 0 $, which is impossible.  The case $ \{ d_1, e_1 \}$ is analogous, solving for $t_1$.

Now consider the case $\{ a_1, e_1 \}$.  (1) reduces to 
$$  d_1 t_2^{\, 2} + 2 c_1 t_1 t_2 + b_1 t_1^{\, 2} = 0 $$
If $d_1=0$ or $b_1 = 0$ we are done. Otherwise, by the quadratic formula,
$$       t_2 =   {-c_1 t_1 \pm  \sqrt{ (c_1^{\, 2} - d_1 b_1)t_1^2 }\over d_1} $$ 
so splitting depends on analysis of the polynomial under the radical.  This is simply  $F(t_1)$ from Lemma 1 and Lemma 2.  We proceed as follows.  From the relations (\ref{eq:e4}) we see that $a_1$ and $e_1$ are each the product of two linear polynomials in the eleven parameters $s_1, \ldots, e$.  Thus $a_1 = 0 = e_1$ leads to four cases, each of which is a system of two linear equations.  This system may be solved, allowing some of the $s_i$ to be replaced with  others.  This greatly simplifies the expressions in (\ref{eq:e4}) for $c_1, d_1$, and $b_1$. Two cases lead to the contradiction $e=0$. In the other two,  we have $F(t_1) = 16 \, e^2 \, s_2^2 \, t_1^2$.  Therefore, $F(t_1)$ is a perfect square and we are done.

The final case $\{ b_1, d_1 \}$ is similar to $\{ a_1, e_1 \}$.  (1) reduces to 
$$  a_1 t_1^{\, 2} t_2^{\, 2}  + 2 c_1 t_1 t_2  + e_1 = 0  $$
If $a_1=0$ we are done. Otherwise, by the quadratic formula,
$$      t_1 \, t_2 =   {-c_1  \pm  \sqrt{ c_1^{\, 2} - a_1 e_1 }\over a_1} $$
$c_1^{\, 2} - a_1 e_1 = F(t_1)/t_1^2$. 
Once again we use the relations (\ref{eq:e4}) to produce four cases.  As before, the solution of two linear systems leads to the contradiction $e=0$; in the other two we have $F(t_1) = 16 \, e^2 \, s_2^2 \, t_1^2$.  Therefore, $c_1^{\, 2} - a_1 e_1$ is a perfect square and we are done.
  $\mqed$ 

\section{List of split cases for $t_2$ a rational function of $t_1$}

The cases and subcases in the previous section may seem bewildering.  We have written a program in a computer algebra system to summarize the details of the four split cases for $t_2$ a rational function of $t_1$.  Recall from Lemma 2 that when this occurs, we have $a_1 b_1=0$ and $d_1 e_1=0$.  This leads to four cases $a_1 = 0$ and $ d_1 = 0$; $a_1 = 0$ and $  e_1 = 0$; $b_1 = 0$ and $ d_1 = 0$; $b_1 = 0$ and $ e_1 = 0$.
 
As in the proof of Lemma 3, cases $\{ a_1, e_1 \}$ and $\{ b_1, d_1 \}$, we use relations (\ref{eq:e4}) to produce systems of two linear equations.  This yields substitutions for one $s_i$ in terms of others, and produces four main cases, each with two subcases. The table shows the resulting $F(t_1)$ and  $t_2$ in terms of $t_1$.

$a_1=0, \, d_1=0:$
\[  \hspace*{8mm}
    \begin{array}{l}
				
  \begin{array}{llll} \hspace*{1mm} s_7  =   \ \ s_5,  \hspace*{2mm} & s_2  =   e\!: \hspace*{3mm} & F(t_1)=16 e^2 s_5^2 t_1^2, \hspace*{3mm} & t_2 =( s_5 t_1^2 + e t_1^2 - s_5 + e  ) / ( 2 s_5 t_1  )\hspace*{1mm}\end{array}      \\
  \begin{array}{llll} \hspace*{1mm} s_7  =  - s_5, \hspace*{2mm} & s_2  =   e\!: \hspace*{3mm} & F(t_1)=16 e^2 s_5^2 t_1^2, \hspace*{3mm} & t_2 =( s_5 t_1^2 - e t_1^2 - s_5 - e  ) / ( 2 s_5 t_1  )\hspace*{1mm}\end{array}   
				
	\end{array}                      
	\]   
$a_1=0,\, e_1=0:$
\[  \hspace*{-6mm}
    \begin{array}{l}
				
  \begin{array}{llll} \hspace*{-16mm} s_7  =   \  s_2,  \hspace*{2mm} & s_5  =   e\!: \hspace*{4mm} & F(t_1)=16 e^2 s_2^2 t_1^2, \hspace*{3mm} & t_2 = t_1 \,{ \textstyle s_2 + e  \over  \textstyle s_2 - e } \hspace*{1mm} \textrm{ or} \end{array}      \\
 \begin{array}{llll} \hspace*{-3mm} \   \hspace*{18mm} & \  \hspace*{18mm} & \  \hspace*{14mm} & t_2 = t_1\hspace*{1mm}\end{array}     \\
  \begin{array}{llll} \hspace*{-16mm} s_7  =  s_2, \hspace*{2mm} & s_5  =  -e\!: \hspace*{2mm} & F(t_1) = 16 e^2 s_2^2 t_1^2, \hspace*{3mm} & t_2 = t_1 \, { \textstyle s_2 + e  \over  \textstyle s_2 - e } \hspace*{1mm} \textrm{ or} \end{array} \\
\begin{array}{llll} \hspace*{-3mm} \   \hspace*{18mm} & \  \hspace*{18mm} & \  \hspace*{14mm} & t_2 =t_1\hspace*{1mm}\end{array}    
				
  \end{array}                      
	\]
 $b_1=0, \, d_1=0:$
\[  \hspace*{11mm}
    \begin{array}{l}
				
  \begin{array}{llll} \hspace*{-32mm} s_7  =   -s_2,  \hspace*{2mm} & s_5 = e\!: \hspace*{3mm} & F(t_1)=16 e^2 s_2^2 t_1^2, \hspace*{3mm}  & t_2 = -1 / t_1 \hspace*{1mm}  \textrm{ or}\end{array}      \\
  \begin{array}{llll} \hspace*{12mm} \ \hspace*{18mm} & \hspace*{6mm} & \hspace*{-1mm} & t_2 = { \textstyle{e + s_2}  \over \textstyle{ t_1 ( e-s_2) } } \hspace*{1mm}\end{array}  \\
 \begin{array}{llll} \hspace*{-32mm} s_7  =   -s_2,  \hspace*{1mm} & s_5 = -e\!: \hspace*{1mm} & F(t_1)=16 e^2 s_2^2 t_1^2, \hspace*{3mm}  & t_2 = -1 / t_1 \hspace*{1mm}  \textrm{ or}\end{array}      \\
  \begin{array}{llll} \hspace*{12mm} \ \hspace*{18mm} & \hspace*{6mm} & \hspace*{-1mm} & t_2 = { \textstyle{e + s_2}  \over \textstyle{ t_1 ( e-s_2) } } \hspace*{1mm}\end{array}   
				
	\end{array}                      
	\]
 $b_1=0, \, e_1=0:$
\[  \hspace*{11mm}
    \begin{array}{l}
				
  \begin{array}{llll} \hspace*{0mm} s_7  =   s_5,  \hspace*{2mm} & s_2 = -e\!: \hspace*{3mm} & F(t_1)=16 e^2 s_5^2 t_1^2, \hspace*{3mm}  & t_2 = 0\hspace*{1mm}  \textrm{ (degenerate) or}\end{array}      \\
  \begin{array}{llll} \hspace*{18mm} \ \hspace*{18mm} & \hspace*{18mm} & \hspace*{13mm} & t_2 = -2 s_5 t_1 / ( s_5 t_1^2 + e t_1^2 - s_5 + e  ) \hspace*{1mm}\end{array}  \\
 \begin{array}{llll} \hspace*{0mm} s_7  =   -s_5,  \hspace*{1mm} & s_2 = -e\!: \hspace*{1mm} & F(t_1) = 16 e^2 s_5^2 t_1^2, \hspace*{3mm}  & t_2 = -2 s_5 t_1 / ( s_5 t_1^2 - e t_1^2 - s_5 - e  ) \hspace*{1mm}  \textrm{ or}\end{array}      \\
  \begin{array}{llll} \hspace*{18mm} \ \hspace*{18mm} & \hspace*{18mm} & \hspace*{13mm} & t_2 = 0 \hspace*{1mm} \textrm{(degenerate)}\end{array}   
				
	\end{array}                      
	\]

Some of the cases above lead to  degenerate solutions, such as $s_7 = s_5, s_2 =e$. This is a ``kite", which was discussed in \cite{adg}.  On the other hand, kites can be part of a non-degenerate configuration if other conditions hold.   Bricard \cite{bric} distinguished two types of (non-degenerate) split solutions.   He was focused on the octahedra.  His {\it Case two} corresponds to two quadrilaterals similar, the third a parallelogram.  {\it Case three}  corresponds to all three  quadrilaterals being parallelograms.  {\it Case one} is non-split, which we now address.
\section{The main theorem}
\mythm{Assuming flexibility, if none of the equations $(\ref{eq:e1}) - (\ref{eq:e3})$ split, then each of the variables $t_i$ is a rational function  of the other two.} 

The main step in the proof of Theorem 1 is the following lemma.  Although considered to be well known, we can find neither proof nor even precise statement of it, so we include it for completeness.
\mylem{ Let $f$ and $g$ be univariate polynomials over some field, say $f = a_n x^n + \ldots +a_0$ and $g = b_m x^m + \ldots + b_0$, where $a_n b_m  a_0 b_0 \ne 0$.  Let $S$ be their Sylvester Resultant matrix, $N \times N$, where $N = n+m$.  If the rank of $S = N-1$, then there exists a polynomial $h(x)$ of degree $1$, whose coefficients are rational functions of $\{ a_i, b_j \}$, satisfied by all the common roots of f and g.}
\myproof  In other words, $x$ is rational function of the coefficients $\{ a_i, b_j \}$.  Note that 0 is not a common root.  We may assume $N \ge 3$.

  We assume familiarity with the basic facts about the Sylvester Resultant.  Since the rank is $N-1$, we may perform row and column operations on $S$ until we have
\[  \hspace*{-4mm}
 S' =  \left[ \begin{array}{l}
				
                    \begin{array}{llllll} \hspace*{2mm} 1 \hspace*{2mm} & 0\hspace*{2mm} & 0\hspace*{2mm} & \ldots\hspace*{2mm} & \hspace*{2mm} 0 \hspace*{2mm} & \hspace*{2mm} c_1 \end{array} \\
		 \begin{array}{llllll} \hspace*{2mm} 0\hspace*{2mm} & 1\hspace*{2mm} & 0\hspace*{2mm} & \ldots\hspace*{2mm} & \hspace*{2mm} 0 \hspace*{2mm} & 
 \hspace*{2mm} c_2  \end{array} \\
		\begin{array}{llllll} \hspace*{2mm} 0\hspace*{2mm} & 0\hspace*{2mm} & 1\hspace*{2mm} & \ldots\hspace*{2mm} & \hspace*{2mm} 0 \hspace*{2mm} & 
\hspace*{2mm} c_3  \end{array} \\
		\begin{array}{llllll} \hspace*{2mm} \ldots\hspace*{2mm} & \ldots\hspace*{2mm} & \ldots\hspace*{2mm} & \ldots\hspace*{2mm} &  &  \hspace*{-1mm} \ldots \hspace*{2mm}  \end{array} \\
		\begin{array}{llllll} \hspace*{2mm} 0 \hspace*{2mm} & 0\hspace*{2mm} & 0\hspace*{2mm} & \ldots\hspace*{2mm} & \hspace*{2mm} 1\hspace*{2mm} & \hspace*{3mm} c_r \end{array} \\
		\begin{array}{llllllll}  \hspace*{2mm} 0 \hspace*{2mm} & 0\hspace*{2mm} & 0\hspace*{2mm} &   \ldots\hspace*{2mm}  &  \hspace*{2mm} 0  \hspace*{2mm} & \hspace*{3mm}0 \end{array}  				
				\end{array}  \right]
	\]
where $r= N-1$.  All of the $c_i $ are rational combinations of the original coefficients $\{ a_i, b_j \}$.  For a common root $x$, the column vector 

$$
\left[ \begin{array}{llllll} \hspace*{2mm} x^{N-1} \hspace*{2mm} & x^{N-2} \hspace*{2mm} & \ldots \hspace*{2mm} & \ x^2 \hspace*{2mm} & x \hspace*{2mm} & 1  \end{array}  \right]^T $$
is in the kernel of the original $S$.   Since column swaps may have been made, the transformed vector 
$$
p = \left[ \begin{array}{llll} \hspace*{2mm}  x^{e_0} \hspace*{2mm} & x^{e_1} \hspace*{2mm} & \ldots \hspace*{2mm} & x^{e_{N-1}}  \end{array}  \right]^T $$ is in the kernel of  $S'$. The exponents are a permutation of the set $\{ 0, 1, 2, \dots, N-1 \}$. If no column swaps were made, the permutation is the identity map and $e_{N-1} = 0$.

If we multiply the matrix $S'$ by the vector $p$ we must get 0.  That produces $N-1$ equations, each a sum of two terms set to 0 (none of the $c_i$ can be 0 as 0 is not a common root). 
We distinguish three cases, according to $e_{N-1} = 0, 1$, or $k > 1$.   In the first case, one of the equations is $x + c_{N-1}=0$, done. In the second case, one of the equations is 
 $1 + c_{N-1} x=0$, done.   In the third case, two of the equations are $x + c_j x^k=0$ and $1 + c_i x^k=0$.  Solve for $x^k$ in one equation, plug into the other, done.
\mqed

{\it Remark:}  This theorem can be generalized to the situation where the rank of the Sylvester matrix is $ < N-1$, but we don't need that here.

To use Lemma 4, we apply the Sylvester resultant method to equations $(\ref{eq:e2})$ and $(\ref{eq:e3})$ to eliminate $t_3$. This Sylvester matrix is $4 \times 4$:
%\smallskip
\addtocounter {equation} {1}
\[  \hspace*{16mm}
  \left[ \begin{array}{l}
		\begin{array}{llll} a_2 t_2^2 + d_2  & \hspace*{12mm} 2 c_2 t_2  & \hspace*{12mm}    b_2 t_2^2 + e_2  & \hspace*{12mm} 0 \end{array} \\
\begin{array}{rrrr} \hspace*{12mm}0  & \hspace*{10mm}  a_2 t_2^2 + d_2  & \hspace*{12mm} 2 c_2 t_2   & \hspace*{10mm}    b_2 t_2^2 + e_2  \end{array} \\
\begin{array}{llll} a_3 t_1^2 + d_3  & \hspace*{12mm}  2 c_3 t_1  & \hspace*{12mm}   b_3 t_1^2 + e_3   & \hspace*{12mm}  0 \end{array} \\
\begin{array}{rrrr} \hspace*{12mm} 0  & \hspace*{10mm}  a_3 t_1^2 + d_3   & \hspace*{12mm}  2 c_3 t_1  & \hspace*{10mm}   b_3 t_1^2 + e_3  \end{array}	 	
	 \end{array}  \right]   \hspace{25mm}  ( \arabic{equation} )
  \]
\mylem{The rank of the Sylvester matrix (7) is 3 almost everywhere.}
\myproof Recall that the $t_i$ are functions of time that are not 0 on any nontrivial interval.  

As equations in $t_3$, $(\ref{eq:e2})$ and $(\ref{eq:e3})$ are quadratic.  The leading coefficients are $a_i t^2 + d_i$ and the ``constant" terms are $b_i t^2 + e_i$ (where $t$ is $t_2$ or $t_1$).  None of these can vanish, as then that equation would split (see Lemma 3).  The hypotheses of Lemma 4 are satisfied.
  
  Since for all values of $t_1$ and $t_2$ in some interval  equations $(\ref{eq:e2})$ and $(\ref{eq:e3})$ have common root(s), the rank is either 0, 1, 2, or 3.  The rank is obviously not 0, as for example $c_2 \ne 0$.  If the rank were 1, then every $2 \times 2$ minor would have determinant 0.  But the upper left $2 \times 2$ minor has determinant $(a_2 t_2^2 + d_2)^2$.  This cannot be 0, as none of the equations split (see Lemma 3). 

If the rank were 2, then every $3 \times 3$ minor would have determinant 0.   Consider then the  minor formed by rows 2, 3, 4, and columns 1, 2, 3.  
\[  \centering
  \left[ \begin{array}{l}
               \begin{array}{lll} \hspace*{12mm}0 &  \hspace*{10mm} a_2 t_2^2 + d_2  & \hspace*{12mm} 2 c_2 t_2    \end{array} \\
               \begin{array}{lll} a_3 t_1^2 + d_3  & \hspace*{12mm}  2 c_3 t_1  & \hspace*{12mm}   b_3 t_1^2 + e_3  \end{array} \\
              \begin{array}{rrr} \hspace*{12mm} 0  & \hspace*{10mm}  a_3 t_1^2 + d_3   & \hspace*{12mm}  2 c_3 t_1  \end{array} \\
         \end{array}  \right] 
  \]
  
Its determinant is $(a_3 t_1^2 + d_3) ( - c_3 a_2 t_1 t_2^2 + c_2 a_3 t_1^2 t_2 + c_2 d_3 t_2 - c_3 d_2 t_1)$.    If this is 0, the second factor must be 0 (by Lemma 3).  Examining the second factor, we distinguish 3 cases:
\begin{itemize}
\item $a_2=0$ and $a_3=0$:   Then again by Lemma 3, $d_2 \ne 0$ and $d_3 \ne 0$.  We immediately solve for $t_2$ as a rational function of $t_1$, contradiction.
\item $a_2 \ne 0$:  In the second factor, solve for $t_2^2$ as a function of $t_2$ (to the first power only) and $t_1$.   Plug this into equation $(\ref{eq:e1})$ to  obtain

%\smallskip
\hspace*{5mm} $  (c_2 a_1 a_3 t_1^4  + c_2 a_1 d_3 t_1^2 + c_2 d_1 a_3 t_1^2 + 2 c_1 c_3 a_2 t_1^2 + c_2 d_1 d_3)t_2 $ 

%\smallskip
\hspace*{56mm} -  $(a_1 d_2 t_1^2 + b_1 a_2 t_1^2 -  d_1 d_2 +  e_1 a_2)c_3 t_1  = 0$

Unless the coefficient of $t_2 = 0$, we can solve for $t_2$ as a rational function of $t_1$, contradiction.  Therefore both expressions in parentheses are 0.  These are both polynomial functions of $t_1$ so their coefficients relative to $t_1$ must be 0.  We immediately see then that $a_1 a_3=0$ (coefficient of $t_1^4$) and $d_1 d_3=0$ (coefficient of $t_1^0$).  Again by Lemma 3, that yields only two possibilities: $a_1=0, d_3=0$ or $a_3=0, d_1=0$.  We are soon led to contradictions, such as $a_2$ must be 0, in both cases. The details are left to the reader.
\item $a_3 \ne 0$:  Exactly like the previous case, only solve for $t_1^2$ instead of $t_2^2$.
\end{itemize} 
This competes the proof that the rank of  Sylvester matrix is 3, except for isolated times when $t_1$ or $t_2$ could be 0. \mqed.
 
 The  proof of Theorem 1 is now  easy:  since the rank of the Sylvester matrix is 3, use Lemma 4 to produce $t_3$ as a rational function of $t_1, t_2$.  By symmetry, any $t_i$ is a rational function of the other two.  \mqed
 
%\medskip
 The proof of Lemma 5 allows us to deduce another result that will soon be of interest:
 \mylem{With the notation of Lemma 5, the rank of the Sylvester matrix (7) is at most 2 (almost everywhere) iff equations $(\ref{eq:e2})$ and $(\ref{eq:e3})$ are multiples of each other, by a nonzero rational function of $t_1$ and $t_2$.}
\myproof The ``if" part is obvious.

Suppose the rank is at most 2.  In the proof of Lemma 5 we used the $3 \times 3$ minor formed by rows 2, 3, 4, and columns 1, 2, 3.  The result was that 
$$   (a_2 t_2^2 + d_2)( 2 c_3 t_1) - (  a_3 t_1^2 + d_3 ) ( 2 c_2 t_2 ) = 0   $$    If we also
consider the minor formed by rows 2, 3, 4 and columns 1, 2, 4, we obtain
$$    ( 2 c_2 t_2 )( b_3 t_1^2 + e_3 ) - ( 2 c_3 t_1) (b_2 t_2^2 + e_2 ) = 0   $$  Thus,  ${\textstyle c_3 t_1 }\over  {\textstyle c_2 t_2}$ times equation $(\ref{eq:e2})$  equals equation $(\ref{eq:e3})$. \mqed
\mycor{If none of the equations $(\ref{eq:e1}) - (\ref{eq:e3})$ split, then for almost all values of $t_1, t_2$, equations $(\ref{eq:e2})$ and $(\ref{eq:e3})$, thought of as equations in $t_3$, do not  have two roots in common.}
\myproof  If they had two roots in common, they would be multiples of each other and the Sylvester rank would be no more than 2, contradicting Lemma 3.  

Of course, the analogous statements can be made about the other pairs of  $(\ref{eq:e1}) - (\ref{eq:e3})$ and the other $t_i$. \mqed

\section{Further analysis of the non-split case, Bricard's case one}

Since none of the equations $(\ref{eq:e1}) - (\ref{eq:e3})$ split, we may use the quadratic formula to solve for, say, $t_2$ and $t_3$ in terms of $t_1$:
%\addtocounter {equation} {1}
\begin{equation} \label{eq:ea}
t_2  =  \frac{-c_1 t_1\pm \sqrt{F(t_1)}}{a_1 t_1^2+d_1}\ 
\end{equation}
%\addtocounter {equation} {1}
\begin{equation} \label{eq:eb}
t_3 = \frac{-c_3 t_1\pm \sqrt{F_1(t_1)}}{a_3 t_1^2+d_3}\ 
\end{equation}

From Theorem 1 we know that $t_3$ is a rational function of $t_1$ and $t_2$.  Therefore

\[ \frac{-c_3 t_1\pm \sqrt{F_1(t_1)}}{a_3 t_1^2+d_3}
= \phi\left(t_1,\frac{-c_1 t_1\pm \sqrt{F(t_1)}}{a_1 t_1^2+d_1}\right) \ 
\]
where $\phi$ denotes a rational function.

Expand $\phi$ and collect terms.  This yields an expression

\[ P\sqrt{F(t_1)} + Q \sqrt{F_1(t_1)} = L  + M \sqrt{F(t_1) F_1(t_1)} \]
and after squaring, eventually  
\[ F(t_1) F_1(t_1) = \frac{R^2}{S^2} \]
for some polynomials $P, Q, L, M, R, S$.
But again, as in Lemma 1, we are over a UFD, so we have proven 
\mythm{With the notation of $(\ref{eq:ea}) - (\ref{eq:eb})$, in the non-split case the product $F(t_1) F_1(t_1)$ is a perfect square, but neither $F(t_1)$ nor $F_1(t_1)$ is a perfect square.} \mqed

Obviously, the same statement is true for the analogous polynomials $F(t_2), F_1(t_2)$, $F(t_3), F_1(t_3).$\footnote{The notation follows Bricard's. Strictly speaking, this is not one function $F$ or $F_1$ being applied to different $t_i$ as the parameters vary.}

  Recall from $(\ref{eq:e6})$ that the $F$ polynomials are in general quartic with no cubic or linear terms:
%\addtocounter {equation} {1}
\begin{equation} \label{eq:e18}
   F(t_1) = - a_1 b_1 t_1^4 + (c_1^2 - a_1 e_1 - b_1 d_1)t_1^2 + d_1 e_1 
\end{equation}
%\addtocounter {equation} {1}
\begin{equation} \label{eq:e19}
F_1(t_1) = - a_3 b_3 t_1^4 + (c_3^2 - a_3 e_3 - b_3 d_3)t_1^2 + d_3 e_3 
\end{equation}  
However, it is possible that, say, $a_1 = 0$, reducing  $F(t_1)$ to a quadratic.  Let us abbreviate $F(t_1) \equiv F,  F_1(t_1) \equiv F_1$.  We distinguish three cases:

\begin{itemize}
\item Both $F$ and $F_1$ are quartic.

\item Both $F$ and $F_1$ are quadratic.

\item One of $F$ and $F_1$ is quartic and one is quadratic.
\end{itemize}
 
 Bricard seems to have missed the possibility of the third case, which we call ``quart-quad."   As he is mostly concerned with octahedra, perhaps he eliminated that case by some three dimensional argument.  
We used {\it Solve} and found no non-split solutions of quart-quad.  Motivated by these experiments, we found a purely algebraic proof of the next Theorem: 
\mythm{If one of $F$ and $F_1$ is quartic and one is quadratic, then we have a split case.}  
\myproof Suppose without loss of generality that $F$ is quadratic, so $a_1 b_1 = 0$.  $F$ and $F_1$ can be factored in some extension field, yielding
\begin{equation}
   F = p(t_1 - \alpha_1)(t_1 - \alpha_2)
\end{equation}
\begin{equation} 
   F_1 = q (t_1 - \beta_1)(t_1 + \beta_1) (t_1 - \beta_2)(t_1 + \beta_2) 
\end{equation}
Since neither $F$ nor $F_1$ is a perfect square but their product is, relabeling if necessary we must have that $\alpha_1 = \beta_1, \alpha_2 = -\beta_1, \beta_2 = 0$.   Thus 
in $F_1$, $d_3 e_3 = 0$.  
Since $F \cdot  F_1$ is a perfect square (in the polynomial ring)  we have $F_1 = s^2 \, t_1^2 \, F$, where $s$ is some polynomial.  Therefore,  
     $\sqrt F_1 = s t_1 \sqrt F$.   
     
 Since $a_1 b_1 = 0$ and  $d_3  e_3 = 0$,  there are four cases.  Let's consider first  $a_1 = 0$ and  $d_3 = 0$.   Then we have  from (\ref{eq:ea}), (\ref{eq:eb})
 \begin{eqnarray}   \label{eq:eg}
       t_2 =    (- c_1  t_1 + \sqrt F) /  d_1 \\    
        t_3 =    (- c_3  t_1 + \sqrt F_1) / ( a_3  t_1^2)
 \end{eqnarray}
 
\noindent
In the second equation, replace $\sqrt F_1$ with $s  \, t_1 \, \sqrt F$.  Simplifying, we have 
      $a_3  t_1  t_3 =  - c_3 + s \sqrt F$.
 Now solve for $\sqrt F$ in the other equation, plug into the above.  We get 
\begin{equation} \label{eq:ee}
    s  \, d_1 \, t_2 = s  \, c_1  \, t_1 +  a_3  t_1  t_3 +  c_3 
\end {equation}

So $t_2$ is a linear function of  $t_1$ or $t_3$.  Now, we know from Theorem 1 that any $t_i$ is a rational function of the others, so (\ref{eq:ee}) may not seem 
surprising.  However, the important fact is that all the exponents are 1.
 
   Solve (\ref{eq:ee}) for $t_2$ and plug that into equation (2).   Collecting terms yields: 
\begin{equation} \label{eq:ef}
    m_1  t_1^2 + m_2  t_1 + m_3 = 0 
\end {equation}
The $m_i$ are polynomials in  $t_3$, up to degree 4.   

Suppose first that no $m_i = 0$.   The key point is that this equation is quadratic in  $t_1$, which
is because (\ref{eq:ee}) is linear in each $t_i$.  Since it is quadratic in  $t_1$ we can apply the same logic used in the proof of Lemma 5, Lemma 6, and Corollary 1 to the pair of 
equations  (\ref{eq:ef}) and  (3).  There must be common root(s).  If the rank of the Sylvester matrix is three, then  $t_1$ is a rational 
function of  $t_3$ (earlier theorem), so we have a split case.  If the rank is less than three, then the two equations are multiples of each other (by a polynomial in the
eleven parameters).  However, this is impossible because the constant (degree zero) terms in the two equations are  
\begin{equation*}
m_3 = s^2  d_1^2  d_2  t_3^2 +  c_3^2 a2  t_3^2 + 2s  c_2  c_3  d_1  t_3 + s^2  d_1^2 e2 +  b_2  c_3^2 , \  \  
 \text{and}  \  \  d_3  t_3^2 +  e_3
 \end{equation*} 
Note that the coefficient of  $t_3$ in the first is $2s  c_2  c_3  d_1$ but there is no  $t_3$ term in the second.  Thus, $s  c_2  c_3  d_1 = 0$, but $ s  c_2  c_3  d_1$ can't be 0 unless  $d_1 = 0$, which implies splitting, since  $a_1 = 0$.   

   Now suppose that some $m_i = 0$.  These are equations in $t_3$ alone, so $m_i = 0$ implies every coefficient in it is 0.  As shown above, this is impossible for $m_3$ (without splitting).  The vanishing of $m_2$ is irrelevant, as  (\ref{eq:ef}) remains quadratic. If only $m_1 = 0$, then $t_1$ is a rational function of $t_3$, hence a splitting. This  proves the first case, that  $a_1 = 0 =  d_3$ implies splitting.
 
 For the other three cases, one  of them is just as above,  but the other two seem harder, because we no longer have the simple monomial denominators of (\ref{eq:eg}).  However, recall that in forming equations $(\ref{eq:e1})-(\ref{eq:e3})$ we may use  cotangent as well as tangent, which means we can replace $t_i$ with $t_i^{-1}$.  That has the effect of switching $a_i \leftrightarrow e_i$ and $b_i \leftrightarrow d_i$, which reinstates the needed monomial denominators.     \mqed

  Suppose now that both $F$ and $F_1$ are quartic.  In a splitting field we have \\

\centerline{ $ F = p (t_1 - \alpha_1)(t_1 + \alpha_1) (t_1 - \alpha_2)(t_1 + \alpha_2)  $ }

\centerline{ $  F_1 = q (t_1 - \beta_1)(t_1 + \beta_1) (t_1 - \beta_2)(t_1 + \beta_2)  $}

But since $F F_1$ is a perfect square, each $\alpha_i$ must equal some $\pm \beta_j$. Therefore, $F$ and $F_1$ are multiples of each other.  The same is true if both $F$ and $F_1$ are quadratic.  This is the key fact in the proof of:
\mythm{Referring to Figure 1 for the angles $\theta_1, \theta_2$ opposite to $\alpha$, in a non-split flexible case we have that $\cos(\theta_1) = \pm \cos(\theta_2)$.  The same cosine relation is true for the angles  $CDG$ and $IHE$ opposite to $\gamma$ (technically $\pi - \gamma$) and for $HIB$ and $AGF$ opposite to $\beta$ (and $\pi - \beta$).}
\myproof We emphasize that this is true throughout the flex.  Unlike Bricard's rather specialized geometric argument, we give an algebraic proof.

Since it is a non-split case we know from above that $F$ and $F_1$ are multiples of each other; let $r_1$ be the ratio. Comparing coefficients, 
 \begin{align}  
d_1 e_1 - r_1 d_3 e_3 = 0 \hspace{3cm} \\
b_1 a_1 - r_1 b_3 a_3 = 0  \hspace{3cm} \\
(c_1^2 - a_1 e_1 - b_1 d_1) - r_1(c_3^2 - a_3 e_3 - b_3 d_3) = 0 \hspace{1cm}
\end{align}
Draw lines $GE$ and $DB$.  $GE$ is on two triangles, one containing $\alpha$, one containing $\theta_1$.  $DB$ is on two triangles, one containing $\alpha$, one containing $\theta_2$.  From the law of cosines we deduce:
\begin{align} 
s_7^2 + e^2 - 2 \, e \, s_7 \cos(\alpha) - (s_5^2 + s_2^2 - 2 \, s_5 \, s_2 \cos(\theta_1)) = 0 \hspace{1cm} \\
s_1^2 + b^2 - 2 \, b \, s_1 \cos(\alpha) - (s_4^2 + s_9^2 - 2 \, s_4 \, s_9 \cos(\theta_2)) = 0 \hspace{1cm}
\end{align}
If we consider the cosines as abstract variables, $(18) - (22)$ is a system of five polynomial equations.  Plug in relations (\ref{eq:e4}) for $a_1, \ldots, d_3$.  Using resultants, we can eliminate any four variables. If we eliminate $s_5,  s_2,  r_1,    \cos(\alpha)$, the resultant is quite simple and has these factors:  
$$\cos(\theta_1),  s_1,  b,  s_7,  s_4, s_9,  e,  \cos^2(\theta_2) - \cos^2(\theta_1) $$
As the resultant must vanish, at least one of these factors must be 0.  $\cos(\theta_1)$ can't be 0, as then $ \cos(\alpha)$ would be a constant.  The only choice is that $ \cos^2(\theta_2) = \cos^2(\theta_1) $.   

The other cases are similar. \mqed

Using Theorem 4 we can form a system of six equations to effectively describe the non-split case.  Assume first that $\cos(\theta_1) = \cos(\theta_2)$.  If we eliminate $\cos(\theta_1)$ from equations (21) and (22), we obtain an expression involving $\cos(\alpha)$ that must be 0, of the form $A  \cos(\alpha) + B$.  As we assume non-degeneracy, this can only be true if $A$ and $B$ are both 0.  We repeat the argument with the two other quadrilaterals, yielding the following six equations (set each to 0):
 \begin{align*}   
e s_4 s_7 s_9 - b s_1 s_2 s_5 ,  \hspace{4.5cm} (23)    \hspace{1.7cm}  \\
- s_2 s_5 s_9^2 - s_4 s_7^2 s_9 + s_4 s_5^2 s_9 + s_2^2 s_4 s_9 - e^2 s_4 s_9 - s_2 s_4^2 s_5 + s_1^2 s_2 s_5 + b^2 s_2 s_5 ,  \hspace{2cm} \\
b s_6 s_8 s_9 + e s_1 s_3 s_4 - b s_1 s_3 s_4 ,  \hspace{6cm}\\
s_6 s_8 s_9^2 + s_1 s_4 s_8^2 - s_4^2 s_6 s_8 - s_1^2 s_6 s_8 + b^2 s_6 s_8 + s_1 s_4 s_6^2 - s_1 s_3^2 s_4 - e^2 s_1 s_4 + 2 b e s_1 s_4   
- b^2 s_1 s_4 ,  \\
e s_5 s_7 s_8 - b s_5 s_7 s_8 - e s_2 s_3 s_6 , \hspace{6cm} \\
s_5 s_7 s_8^2 + s_3 s_6 s_7^2 - s_5 s_6^2 s_7 - s_3^2 s_5 s_7 + e^2 s_5 s_7 - 2 b e s_5 s_7 + b^2 s_5 s_7 + s_3 s_5^2 s_6 - s_2^2 s_3 s_6   
- e^2 s_3 s_6 
\end{align*}   
Minor variations result by using $\cos(\theta_1) = -\cos(\theta_2)$, etc.
\section{Flexibility analysis with symbolic software}
The program {\it Solve} was described in \cite{adg} and \cite{fox}.  Here is a brief description.

Let $res$ be the resultant of a system of equations defining a structure, such as $(1)-(3)$.  $res$ is a polynomial in one of the angles, say $t$, and the fifteen parameters $a_1, b_1, \ldots, e_3$, or alternatively, in the eleven side parameters $s_1, s_2, \ldots, e$.  If the structure is flexible, then infinitely many values of $t$ satisfy the polynomial.  The only way this can be is if every coefficient of $t^k$ vanishes.  {\it Solve} examines these coefficients finding ways to kill them one-by-one, usually starting at the top coefficient.  Whenever a way is found to kill the coefficient of $t^k$, that substitution is put on a stack and applied to $res$, creating a polynomial $res'$ of fewer terms and one fewer parameter.  Then {\it Solve} calls itself on the new polynomial $res'$.  This is essentially an enormous tree search.  Many heuristics and techniques are used to keep the search manageable yet effective.

The output of the algorithm is a list of {\it tables} consisting of substitutions of the form $s_i = p(s_1, s_2, \ldots, \hat{s_i}, \ldots)$, where $p \in$ \GF.   
Here is a simple example.  If $res$ were 
$(s_9 s_8 - s_7 s_6)t^2 + ({s_4}^2 - {s_3}^2)t + s_8 - s_6$, one 
solution would be the table of the three relations $s_9 = s_7, s_8 = s_6, s_4 = s_3$. 

 The relations may be described as follows:
Partition  the set of $N$ parameters into nonempty subsets $X = \{ x_i \}_{i=1}^n$,
$Y = \{ y_j \}_{j=1}^m$, $n+m=N$.  Each relation is an equation $y_j = g(x_{i_1}, x_{i_2}, \ldots)$ where $g$ is a
rational function. A collection of $m$ of these for $j = 1, \ldots, m$ is a {\it solution table} if $res$ evaluated at
them all is 0.  In the example above
$X = \{ s_3, s_6, s_7 \}$ and $Y = \{ s_4, s_8, s_9 \}$.

%\smallskip
\par\noindent {\bf Problem 3: Can all flexible cases be represented by a table of relations in the above sense?}
% \smallskip 

To apply this to the quadrilaterals of Figure 1, we eliminate two of the three angles in equations $(1)-(3)$.  In terms of the eleven side parameters, $res$ has 190981 terms\footnote{By late 2013, this computation takes {\it Fermat} 1.86 minutes on a Mac Mini.}.  From  2006 to 2011, {\it Solve(res)} found many flexible cases of Bricard's types two and three, and many degenerate cases \cite{adg}, \cite{fox}.  Improvements by 2012 yielded the first non-split cases, Bricard's case one.  These were all what we call {\it isohexagons}.   Here is an actual table as computed by {\it Solve}:
\begin{align*}  
        s_9 = ( b \, s_3  ) / ( b - e  )  \hspace{.2cm}  \\
            s_8 =  s_2 ( e  - b   ) / e  \hspace{.6cm}  \\
            s_7 = ( e \, s_1  ) / b   \hspace{1.2cm} \\
       s_6 =  s_1  ( b - e  ) / b  \hspace{.6cm} \\
            s_5 = ( e \, s_3  ) / ( b - e  )  \hspace{.2cm} \\
            s_4 = ( b \, s_2  ) / e    \hspace{1.2cm}  
\end{align*}
Note that $s_8$ is negative; that just means that point $H$ is below the $x$-axis.   They are called isohexagons because a hexagon with equal opposite sides appears in the midst of the flex.  Here are two images (Figure 2), taken from a model made by plugging in numbers for the sides.  The hexagon is outlined in blue.  Note that it does not consist entirely of sides $s_i$.

These isohexagons exhibit a curious kind of symmetry or ``quasi-similarity."  Let $r = (b-e)/e$.

\begin{figure} 
\centering
\includegraphics[scale=0.30]{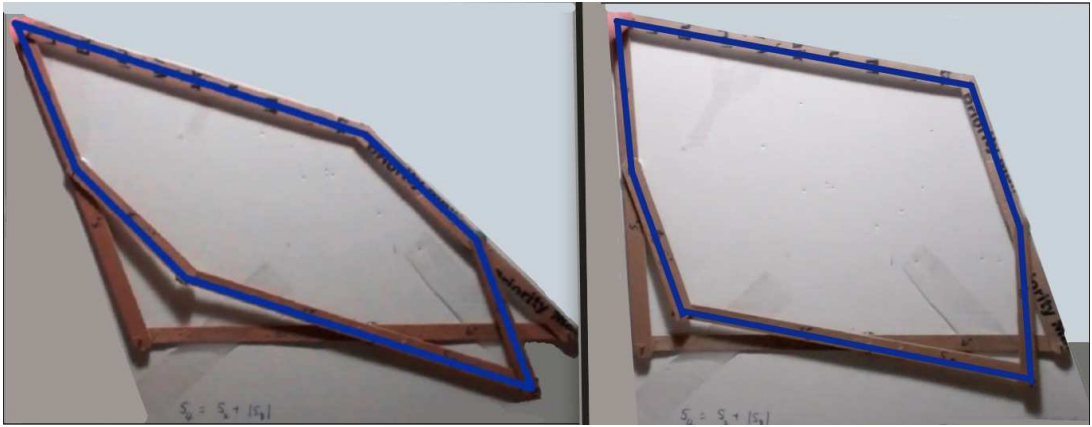}

\caption[ ]{Flexing of case one quadrilaterals, forming ``isohexagon."} 
\end{figure}
Then the table above is equivalent to this:
\[  \hspace*{-5mm}
  \begin{array}{l}
	          \begin{array}{lll} \hspace*{2mm} lower \ left &  \hspace*{4mm} lower \ right &   \hspace*{9mm}   large \end{array} \\
                   \begin{array}{lll} \hspace*{8mm} e  &  \hspace*{15mm} b-e = r \, e &   \hspace*{8mm}   b  = (1+r) \, e \end{array} \\
				\begin{array}{lll} \hspace*{8mm} s_2  &  \hspace*{13mm} s_8  = -r \, s_2 &   \hspace*{8mm} s_4 = (1+r) \, s_2\end{array}  \\
				\begin{array}{lll} \hspace*{8mm} s_5  &  \hspace*{13mm}  s_3  = r \, s_5 &  \hspace*{11mm} s_9 =  (1+r) \, s_5 \end{array}  \\
				\begin{array}{lll} \hspace*{8mm} s_7  &  \hspace*{13mm} s_6  = r \, s_7 &  \hspace*{11mm} s_1 = (1+r) \, s_7 \end{array}  \\
\end{array}				   
	\hspace*{3mm}	     \]
Each column lists the four sides of one of the quadrilaterals.   Note that the sides of the lower right and large quadrilaterals are multiples of those of the lower left, but in an odd shifting pattern; none of these quadrilaterals is similar to another.  From this and other examples we see that the following {\it isohexagon property} holds: for every one of the twelve sides, say $x$, there is a side on a different quadrilateral, say $y$, such that $x/y$ or $y/x$ equals $r, 1+r, -r$, or $-(1+r)$.    Bricard, who studied only octahedra, remarks that in case one  there is an odd symmetry also.  

Thus, the isohexagon may be thought of as the three-quadrilaterals analogue of the octahedra case one.  
We therefore conjectured that all case one examples were isohexagons.  Surprisingly, this is false.   To see why, recall the six equations (23).  All six of these must be 0.  Form a single polynomial 
\addtocounter{equation}{1} 
\begin{equation}  \label{eq:f}
 f = c_5 t^5 + c_4 t^4 + c_3 t^3 + c_2 t^2 + c_1 t + c_0
\end{equation}
where $\{ \, c_i \, \}$ are the six equations in (23) and $t$ is an abstract variable.  Execute {\it Solve(f)}, adding code to suppress split cases.  This finishes very quickly with 136 tables, several of which are striking, such as
\begin{align*}  
s_9 = s_3 ( - b \, e \, s_2^2  + b^2 \, s_2^2  + b \, e^3  - b^2 \, e^2 
					) / ( e^2 \, s_6^2 - e^2 \,  s_3^2  ) \hspace{.1cm}  \\
	 s_8 =  s_2 ( e  - b   ) / e  \hspace{4.2cm}  \\
	 s_7 = ( e \,  s_6  ) / ( b - e  ) \hspace{3.8cm}  \\
	 s_5 = ( e \,  s_3  ) / ( b - e  ) \hspace{3.8cm}  \\
	 s_4 = ( - b \,  s_2  ) / e  \hspace{4.5cm}  \\
	 s_1 =  s_6 ( b \, e \, s_2^2  - b^2 \, s_2^2  - b \, e^3  + b^2 \, e^2 
				) / ( e^2 \, s_6^2 - e^2 \, s_3^2  )   \hspace{.4cm} 
\end{align*}

This is not an isohexagon, as  it does not satisfy the isohexagon property.  Substituting numerical values, we created a model of this case.  Figure 3 shows two snapshots during the flex.

\vspace*{-1mm}
\begin{figure} 
\centering
\includegraphics[scale=0.32]{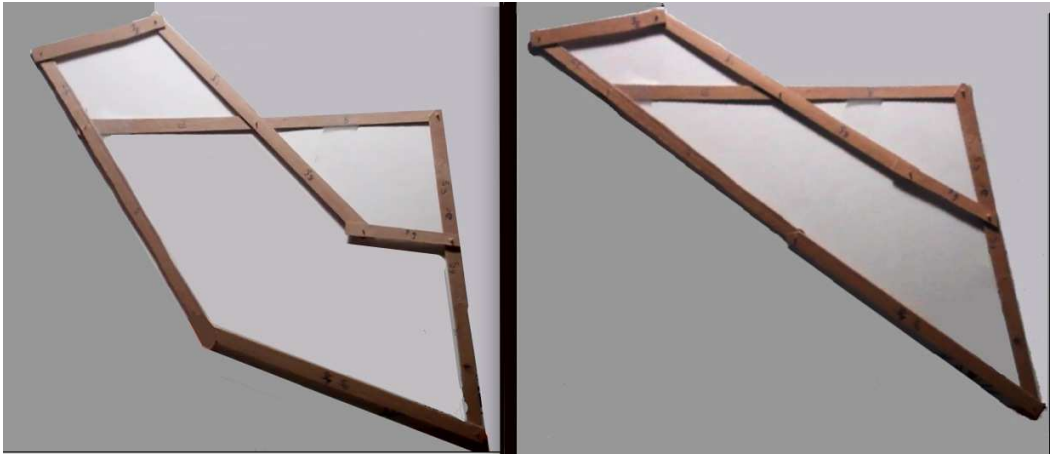}
\caption[ ]{Flexing of case one quadrilaterals, not an isohexagon.} 
\end{figure}
\vspace*{-1.5mm} 
\noindent 
The last case above seems rather complicated due to the $s_1$ and $s_9$ equations.  Notice that  $s_9/s_1 = s_3/s_6$.  To experiment, we removed the $s_9$ and $s_1$ equations, added  $s_9 = s_1 s_3 / s_6$, and plugged the resulting table into the resultant $res$ of 190981 terms.  It did not vanish, but left a polynomial $res'$ of 8803 terms.  Further analysis of this polynomial revealed a surprising ``irrational case":

\vspace*{-7.5mm} 
\begin{align*}  
s_9 =   ( - b \, s_3 ) / ( e - b  ),  \hspace{3cm}  \\
s_8 =  s_2 ( e - b   ) / e ,  \hspace{3.69cm}  \\
s_7 =  e \, s_6   / ( e - b  ),   \hspace{3.62cm}  \\
s_5 =  e \, s_3   / ( e - b  ),  \hspace{3.62cm}  \\
s_4 =  b \, s_2   / e  \hspace{4.83cm}  \\
s_1 = ( - b \, s_6  ) / ( e - b  )  \hspace{3.15cm}  \\
s_6^2 = ( e^2  \, s_3^2 + e^2 \,  s_2^2 - 2 \,  b \,  e  \, s_2^2 + b^2 \,  s_2^2 - e^4 + 2 \,  b  \, e^3 - b^2 \,  e^2  ) / e^2 \hspace{.15cm} 
\end{align*}  (Use the last equation to replace $s_6$ in the earlier relations.)  
This arises because $s_6$ occurs with only even degree in $res'$.  
In the definition of table, the sets $X$ and $Y$ are
\begin{align*}
X = \{ \,  b, e, s_2, s_3 \, \} \hspace{2cm} \\
Y = \{ \,  s_1, s_4, s_5, s_6, s_7, s_8, s_9 \,   \}
\end{align*}
Figure 4 is an image of an instantiation of this, plugging in real numbers for the $X$ (and then $Y$) parameters.

\vspace*{-2mm}
\begin{figure} 
\centering
\includegraphics[scale=0.32]{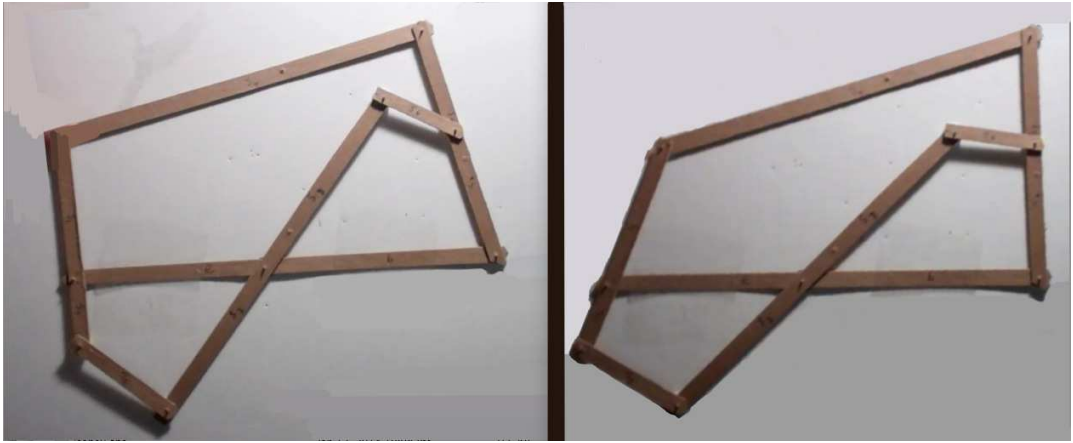}
\caption[ ]{Flexing of case one quadrilaterals, irrational relationship.} 
\end{figure}
\vspace*{-1mm} 
\noindent 

%\bigskip

The polynomial in the definition of $s_6^2$ is easily seen to be not a perfect square.  Therefore, the answer to Problem 3 is ``no."  

The three structures in Figures $2 - 4$ were discovered first using {\it Solve}($f$) for $f$ defined in (\ref{eq:f}).  However, all  three now show up with the latest version of {\it Solve(res)}. This is because when {\it Solve} encounters a polynomial (like $res'$) in which a variable (here $s_6$) occurs with every exponent a multiple of $n$ (here 2), the exponents are divided by $n$ and the algorithm continues.

\section{Conclusion}

This problem of the flexible planar linkages was posed by Bricard in his memoir on the flexible octahedron. He seemed to imply that the two problems would have completely analogous resolutions given that they are described by systems of equations of identical form.  As Bricard pointed out in his memoir, ``it ought to be possible to analyze these equations by purely algebraic means, however the amount of computation required would be daunting". He proceeded therefore to analyze his equations geometrically, arriving at his well known three classes of flexible octahedra. Here, with the help of computer algebra, which a hundred years after Bricard is now a mature field, we were able to carry out this ``daunting" task for the planar mechanism case, and were rewarded by a surprising divergence from Bricard's conclusions.  Although the separation to three classes according to the type of splitting is identical for the two problems, the underlying geometric differences led to 
unexpectedly rich properties for the structures of case one, the case of no splitting, with no analogs in the octahedron. The other two split cases seem to be completely analogous for the two problems.

Earlier we defined three problems:

\noindent
{\it Problem 1:  Find conditions  on the sides under which the quadrilateral arrangement becomes flexible.}
This has been solved.

\noindent
{\it Problem 1$'$:  Find all conditions on the sides under which the quadrilateral arrangement becomes non-degenerate flexible.}  
We do not have a mathematical proof that our list is complete. However we have found cases analogous to all of Bricard's cases for the articulated octahedron, and discovered unexpectedly rich properties for case one, where our algebraic analysis led to two quite different types of flexible structures with no apparent analogy to Bricard's three-dimensional  results.

\noindent
{\it Problem 2:  When is one of these variables, $t_2$, say, a rational function of another $t_j$, or a rational function of both of the other ones $t_1, t_3$? }
This has been solved.  Our proofs are new and algebraic.

\noindent
{\it Problem 3:  Can all flexible cases be represented by a table of relations in [our] sense?}  
No.  However, we believe that our algorithm needs to be modified only by changing the 
 definition of table to allow relations of the form  $s_i^n = p(s_1, s_2, \ldots,$ $ \hat{s_i}, $ $\ldots)$.

 The great success we have had on this project bodes well for future work with more complex structures, as equations describing those structures are also quadratic, based on distances and angles.

 Work is ongoing applying these methods directly to the octahedra and to the cyclo-octane molecule.

%\end{large}

%
% ---- Bibliography ----

\end{document}